\magnification=\magstep1
\tolerance=2000

\baselineskip=17pt

\font\tenopen=msbm10
\font\sevenopen=msbm7
\font\fiveopen=msbm5
\newfam\openfam
\def\openo{\fam\openfam\tenopen}
\textfont\openfam=\tenopen
\scriptfont\openfam=\sevenopen
\scriptscriptfont\openfam=\fiveopen

\font\tenujsym=msam10
\font\sevenujsym=msam7
\font\fiveujsym=msam5
\newfam\ujsymfam

\textfont\ujsymfam=\tenujsym
\scriptfont\ujsymfam=\sevenujsym
\scriptscriptfont\ujsymfam=\fiveujsym

\def\C{{\openo C}}
\def\R{{\openo R}}

\def\kocka{\kern+0.9pt{\sqcup\kern-6.5pt\sqcap}\kern+0.9pt}

\centerline{\bf Addendum to my paper} 

\centerline{\bf ``The Lebesgue summability of trigonometric integrals"} 

\centerline{\bf (J. Math. Anal. Appl. 390 (2012), 188-196)} 

\bigskip

\centerline{\bf Ferenc M\'oricz}

\vglue1cm

We observed that statement (2.9) in Theorem 2 remains valid if condition (2.7) is replaced 
by the weaker condition that 
$$f(t) \in L^1(-T, T) \quad {\rm for\ all} \quad T>0.\eqno(2.7')$$
Furthermore, Theorem 3 also remains valid if condition (2.7) is replaced by (2.7$'$). 

To be more precise, the following Theorems 2$'$ and 3$'$ can be proved in the same way as 
Theorems 2 and 3 are proved, while using Lemmas 2$'$ and 3$'$ below instead of Lemmas 2 and 3. 

\bigskip

\noindent {\bf Theorem 2$'$.} {\it If $f: \R\to \C$ is such that conditions} (2.7$'$) 
{\it and 
$$\lim_{T\to \infty} {1\over T} \int_{|t|<T} |tf(t)| dt =0\eqno(2.8)$$
are satisfied, then we have uniformly in $x\in \R$ that }
$$\lim_{h\downarrow 0} \Big\{{\Delta {\cal L} (x; h)\over 2h} - I_{1/h} (x)\Big\} = 0.\eqno(2.9)$$

\bigskip

\noindent {\bf Theorem 3$'$.} {\it Suppose $f: \R\to \C$ is such that conditions} 
(2.7$'$) {\it and} 
$${1\over T} \int_{|t|<T} |tf(t)| dt \le B \quad for\ all\quad T>T_1, \eqno(2.10)$$
{\it are satisfied, where $B$ and $T_1$ are constants. If the finite limit 
$$\lim_{T\to \infty} I_T (x) : = \lim_{T\to \infty} \int_{|t|<T} f(t) e^{itx} dt = \ell\eqno(2.3)$$
exists at some point $x\in \R$, then} (2.9) {\it holds at this} $x$. 

We emphasize that in this addendum, the definition 
$$\int_\R f(t) {e^{itx}\over it} dx =: {\cal L} (x), \quad x\in \R,\eqno(2.4)$$
is interpreted only formally; that is, the integral in (2.4) may not exist in 
Lebesgue's sense. However, under the conditions in Theorems 2$'$ and 3$'$, the integral in the 
representation 
$${\Delta {\cal L} (x;h)\over 2h} : = \int_\R f(t) e^{itx} {\sin th\over th} dt, \quad h>0,\eqno(2.6)$$
does exist in Lebesgue's sense. 

As we have mentioned above, the proofs of Theorems 2$'$ and 3$'$ hinge on the following Lemmas 2$'$ and 3$'$. 
We note that we essentially use only Part (i) 
in our earlier Lemmas 2 and 3, while we substitute condition (2.7$'$) for (2.7). 

\bigskip

\noindent {\bf Lemma 2$'$.} {\it If $f: \R\to \C$ is such that condition} (2.7$'$) {\it and} (2.8) 
{\it are satisfied, then} 
$$\lim_{T\to \infty} T \int_{|t|> T} \Big|{f(t)\over t}\Big| dt =0.\eqno(3.2)$$

\bigskip

\noindent {\bf Lemma 3$'$.} {\it If $f: \R\to \C$ is such that conditions} (2.7$'$) {\it and} (2.10) {\it are satisfied, then} 
$$T\int_{|t|>T} \Big|{f(t)\over t}\Big| dt \le 4B \quad for\ all\quad T> T_1. \eqno(3.8)$$

We note that the converse implications (3.2) $\Rightarrow$ (2.8) in Lemma 2$'$ and (3.8) $\Rightarrow$ 
(2.10) in Lemma 3$'$ do hold under the supplementary condition (2.7). But we do not need these converse implications in 
the proofs of Theorems 2$'$ and 3$'$. 

\bye